\documentclass{amsart}
\usepackage{amssymb,latexsym,amsmath,amscd,amsthm,amsfonts}
\begin{document}

\title [quantum group actions] {Compact quantum group actions on C*-algebras\\ and invariant derivations}
\author{R. Dumitru}
\address{Department of Mathematics, 
University of Cincinnati, Cincinnati, OH and Institute of Mathematics of the Romanian Academy Bucharest, Romania}
\email{dumitrra@email.uc.edu}
\author{C. Peligrad}
\address{Department of Mathematics, 
University of Cincinnati, Cincinnati, OH and Institute of Mathematics of the Romanian Academy Bucharest, Romania}
\email{costel.peligrad@UC.Edu}
\thanks{Research supported by the Taft Foundation. The first author was supported by a Taft Advanced Graduate Fellowship and by a Taft Graduate Enrichment Award. The second author was supported by a Taft Research Travel Grant.}
\maketitle

\begin{abstract} We define the notion of invariant derivation of a C*-algebra under a compact quantum group action and prove that in certain conditions, such derivations are generators of one parameter automorphism groups.
\end{abstract}

\

\

\

\section{Introduction}

\ In this paper we consider actions of a compact quantum group $G=(A,\Delta)$ on a C*-algebra $B$.

\ In Section 2 we collect some preliminary results about the spectral subspaces of such actions.

\ In Section 3 we define the notion of invariant derivation under a compact quantum group action and prove that such derivations are generators of one-parameter automorphism groups if their domain contains all the spectral subspaces (Theorem 3.8) or, if $B=\mathcal{K}(H)$, if the domain contains the fixed point algebra $B^{\delta}$.  

\ Our results extend and improve on results obtained in \cite{bra}, \cite{gdm}, \cite{plgr} for the case of groups.

\

\section{Compact Quantum group actions on C*-algebras}

\ Let $G = (A,\Delta)$ be a compact quantum group, i.e. a unital C*-algebra $A$ and $\Delta : A \to A\otimes A$ a $\ast$ - homomorphism such that

\ $(i)$ $(\Delta \otimes \iota )\Delta = (\iota \otimes \Delta )\Delta$, where $\iota$ is the identity map\\and

\ $(ii)$ $\overline {\Delta(A)(1\otimes A)}=A\otimes A$ and $\overline {\Delta(A)(A\otimes 1)}=A\otimes A$.

\ Let $\widehat G$ denote the dual of $G$, i.e. the set of all unitary equivalence classes of irreducible representations of $G$ (\cite{wor1}, \cite{wor2}).
\ For each $\alpha \in \widehat G$, denote by $\ u ^{\alpha}$ a representative of each class. Let $ \{u_{ij}^{\alpha}\}_{i,j=\overline{1,d_\alpha}}\subset A$ be the matrix elements of $\ u ^{\alpha}$.

\ Set $ \chi_\alpha =\underset{i=1}{\overset{d_{\alpha}}{\sum}} u_{ii}^{\alpha}$. Then $\chi_\alpha$ is called the character of $\alpha$. By (\cite{wor1}, \cite{wor2}) there is a unique invertible operator $F_{\alpha} \in B(H_{\alpha})$, where $H_{\alpha}$ is the finite dimensional Hilbert space of the representation $u^{\alpha}$, that intertwines $u^{\alpha}$ with its double contragradient representation $(u^{\alpha})^{cc}$ such that $tr(F_{\alpha})=tr(F_{\alpha}^{-1})$. Set $M_{\alpha}=tr(F_{\alpha})$.

\ $F_{\alpha}$ can be represented as a matrix $F_{\alpha}=[f_{1}(u_{ij}^{\alpha})]$ where $f_{1}$ is a linear functional on the dense $\ast$-subalgebra $\mathcal{A}\subset A$ spanned (linearly) by $ \{u_{ij}^{\alpha}\}_{\alpha \in \widehat G,i,j}$.

\ If $a\in A$ (or $a\in \mathcal{A}$) and $\xi$ is a linear functional (on $A$ or $\mathcal{A}$), we denote (\cite{wor1}, \cite{wor2})
\[
\ a \ast \xi = (\xi \otimes \iota)(\Delta(a))
\]
\[
\xi \ast a = (\iota \otimes \xi)(\Delta(a))
\]

\ Both $a \ast \xi$ and $\xi \ast a$ are elements of $A$ (respectively $\mathcal{A}$).
 
\ If $h$ is the Haar measure on $A$ (\cite{wor1}), set 
\[
\ h_{\alpha} = M_{\alpha} h\cdot (\chi_{\alpha}\ast f_{1})^{\ast},
\]
where $(h\cdot a)(b)=h(ab)$.

\ Let now $B$ be a C*-algebra and $\delta: B\rightarrow B\otimes A$ be a $\ast$- homomorphism of $B$ into the multiplier algebra of the minimal tensor product $B\otimes A$. If:

\ $(a)$ $\;(\iota \otimes \Delta)\delta = (\delta \otimes \iota)\delta$ and
 
\ $(b)$ $\; \overline {\delta(B)(1\otimes A)}=B\otimes A $\\
then $\delta$ is called an action of $G$ on $B$, or a coaction of $A$ on $B$.

\swapnumbers
\newtheorem{definition}{Definition}[section]
\newtheorem{remark}[definition]{Remark}

\begin{remark} \label{rem} The condition b) above implies that $\delta \in Mor(B,B\otimes A)$ in the sense of (\cite{wor2}, Introduction). Then, by the discussion in \cite{wor2}, $\delta$ can be uniquely extended to a $\ast$-homomorphism $\delta ^{\prime \prime} :\mathcal{M}(B)\to \mathcal{M}(B\otimes A)$. 
\end{remark}

\ Let $P_{\alpha} : B \rightarrow B$ be the linear map $P_{\alpha}(x)=(\iota \otimes h_{\alpha})(\delta (x))$.

\ In particular, if $\alpha =\iota$ (the trivial one-dimensional representation),\\ $P_{\iota}=(\iota \otimes h)\circ \delta$ is a conditional expectation from $B$ to $B^{\delta}=\{ x\in B\mid \delta(x)=x\otimes 1_{A}\}$ which will be called the fixed point algebra of the action.

\begin{definition}\hfill
\begin{enumerate}
\item For every $\alpha \in \widehat{G}$, denote 
\[ B_{\alpha}=\{ P_{\alpha}(x)\mid x\in B\}. \]
\ The subspace $B_{\alpha}\subseteq B$ is called the spectral subspace corresponding to $\alpha$.
\item For $\alpha \in \widehat{G}$, $i, j = 1,...d_{\alpha}$ we denote 
\[ C_{ij}^{\alpha} = M_{\alpha} (u_{ij}^{\alpha} \ast f_{1})^{\ast} \]
\ and
\[ P_{ij}^{\alpha}(x)=(id\otimes h\cdot c_{ij}^{\alpha})(\delta (x)),\; x\in B.\]
\end{enumerate}
\end{definition}

\

\ We collect some properties of these objects in the following Lemma (see \cite{boca}, \cite{podl}, or easy calculations).
\newtheorem{lemma}[definition]{Lemma}
\begin{lemma} \label{L:P} \hfill
\begin{list}{}
\item(i) $P_{ij}^{\alpha}P_{kl}^{\beta}=\delta_{il}\delta_{\alpha \beta}P_{kj}^{\alpha}$, where $\delta_{il}$, $\delta_{\alpha \beta}$ denote the Kronecker symbols.

\item(ii) $P^{\alpha}P_{ij}^{\alpha}=P_{ij}^{\alpha}$.

\item(iii) If $x\in B_{\alpha}$, then $x= \pmb{\sum}_{i=1}^{d_{\alpha}}P_{ii}^{\alpha}(x) $.

\item(iv) The algebraic direct sum $\pmb{\sum}_{\alpha \in \widehat G}B_{\alpha}$  is a dense $\ast$-subalgebra of $B$.

\item(v) $\delta(P_{ij}^{\alpha}(x))=\pmb{\sum}_{l=1}^{d_\alpha}P_{il}^{\alpha}(x)\otimes u_{lj}^{\alpha}$, $\forall x\in B$.
\end{list}
\end{lemma}

\

\ For $\alpha \in \widehat{G}$, we can consider a matricial spectral subspace:
\[ B_{2}(\alpha)=\{ [P_{ij}^{\alpha}(x)]\mid x\in B\}.\]

\ Then $B_{2}(\alpha)\subset B\otimes B(H_{\alpha})$.

\ By Lemma \ref{L:P}, $(iii)$ the map $x\mapsto [P_{ij}^{\alpha}(x)]$ is a linear isomorphism between $B_{\alpha}$ and $B_{2}(\alpha)$.

\begin{definition}
\ If $\delta :B\rightarrow \mathcal{M}(B\otimes A)$ and $\theta: C\rightarrow \mathcal{M}(C\otimes A)$ are actions of $G$ on $B$ and $C$, denote by $\delta \textsf{\textcircled{\scriptsize{T}}}\theta$ the tensor product of the two actions, i.e.
\[(\delta \textsf{\textcircled{\scriptsize{T}}} \theta)(b\otimes c)= \delta(b)_{13} \theta(c)_{23}.\]
\end{definition} 

\

\ We used the leg-numbering notation (see for instance \cite{baaj}). Straightforward calculations show that this is a coaction of $A$ on $B\otimes C$.

\ In particular, if $\iota$ denotes the trivial action of $G$ on $B(H_{\alpha})$, $\iota(x)= x\otimes 1$, $\forall x\in B(H_{\alpha})$, then $\tilde{\delta}=\delta \otimes \iota$ is an action of $G$ on $B\otimes B(H_{\alpha})$.

\begin{remark} \label{R:tilda}
\ For $x\in B_{\alpha}$, let $X=[P_{ij}(x)]\in B_{2}(\alpha)$. Then:
\[ \tilde{\delta}(X)=(X\otimes 1_{A})(1_{B}\otimes u^{\alpha})\]
\end{remark}

\ This follows from Lemma \ref{L:P} $(v)$.

\begin{lemma} \label{L:zero} $B^{\delta}$ is a C*-subalgebra of $B$ and $B^{\delta}\ne (0)$.
\end{lemma}
\begin{proof} It is immediate that $B^{\delta}$is a C*-subalgebra.

\ We prove next that $B^{\delta}\ne (0)$. By Lemma \ref{L:P} $(iv)$ there is $\alpha\in\widehat{G}$ such that $B_{\alpha}\ne (0)$. Let $x\in B_{\alpha}$, $x\ne 0$. Let \[ X=[X_{ij}]\in B_{2}(\alpha)\] as in Remark \ref{R:tilda}. Then $X\ne 0$. By Remark \ref{R:tilda}, we have \[ \tilde{\delta}(X)=(X\otimes 1_{A})(1\otimes u^{\alpha})\].

\ Hence \[ \tilde{\delta}(XX^{\ast})=XX^{\ast}\otimes 1_{A}.\]

\ Therefore each entry of $XX^{\ast}$ belongs to $B^{\delta}$. Since $X\ne 0$, it follows that $B^{\delta}\ne (0)$. 
\end{proof}

\begin{lemma} \label{L:app} If $(e_{\lambda})_{\lambda}$ is an approximate identity of $B^{\delta}$, then $(e_{\lambda})_{\lambda}$ is an approximate identity of $B$.
\end{lemma}
\begin{proof} Let $x\in B_{\alpha}$, for some $\alpha\in \widehat{G}$ and $X=[X_{ij}]\in B_{2}(\alpha)$ as in Remark \ref{R:tilda}. Then, as in the proof of the previous lemma, $XX^{\ast}\in (B\otimes B(H_{\alpha}))^{\tilde{\delta}}$.

\ For each $\lambda$, let $E_{\lambda}=[e_{ij}^{\lambda}]$ with $e_{ii}^{\lambda}=e_{\lambda}$ for every $i$ and $e_{ii}^{\lambda}=0$ if $i\ne j$.

\ Since $(e_{\lambda})_{\lambda}$ is an approximate identity of $B^{\delta}$ and $XX^{\ast}\in (B\otimes B(H_{\alpha}))^{\tilde{\delta}}$, we obviously have: 
\begin{align} \lim_{\lambda}\|E_{\lambda}X-X\|&=\lim_{\lambda}\|(E_{\lambda}X-X)(E_{\lambda}X-X)^{\ast}\|^{\frac{1}{2}}\notag \\&=\lim_{\lambda}\|E_{\lambda}XX^{\ast}E_{\lambda}-E_{\lambda}XX^{\ast}-XX^{\ast}E_{\lambda}+XX^{\ast}\|^{\frac{1}{2}}=0.\notag
\end{align}

\ Since, by Lemma \ref{L:P} $(iii)$, $x=\sum X_{ii}$, it follows that $\lim_{\lambda}\|e_{\lambda}x-x\|=0$. 

\ As $\overline{\underset{\alpha \in \widehat G}{\sum}B_{\alpha}}=B$ (Lemma \ref{L:P} $(iv)$), the proof is complete.
\end{proof}

\ Since $B^{\delta}\ne (0)$, then the minimum dimension that $B^{\delta}$ can have is 1.
\begin{definition} The action $\delta$ is called ergodic if $\dim B^{\delta}=1$.
\end{definition}

\begin{remark} If $\delta$ is ergodic, then $B$ has a unit. Indeed, by Lemma \ref{L:zero}, $B^{\delta}$ is a one-dimensional C*-subalgebra of $B$ and thus it has a unit. By Lemma \ref{L:app} this is a unit of $B$ also.
\end{remark}

\ Interesting examples and results on ergodic actions can be found in \cite{vaes}, \cite{boca}, \cite{wang}.

\

\section{Unbounded derivations commuting with quantum group actions}

\ Let $d:D(d)\rightarrow B$ be a closed, densely defined $\ast$-derivation. Here $D(d)$ is a dense $\ast$-subalgebra of $B$.

\ As it is easy to show, $d \otimes i_{A} : D(d) \otimes A \rightarrow B\otimes A$ is a closable linear operator. Denote by $\overline{d\otimes i_{A}}$ its closure.

\begin{definition} We say that $d$ is $\delta$-invariant (or that $d$ commutes with $\delta$) if for every $x\in D(d)$, it follows that:
\begin{align}
a)\;&\delta (x)\in D(\overline{d\otimes i_{A}})\notag \\
b)\;&\delta(d(x))=(\overline{d\otimes i_{A}})(\delta(x)),\; \forall x\in D(d)  \notag
\end{align}
\end{definition}

In order to give an equivalent formulation of the concept of invariant derivation, we need the following:
\begin{definition}If $\delta_{1}:B_{1}\rightarrow \mathcal{M}(B_{1}\otimes A)$ and $\delta_{2}:B_{2}\rightarrow \mathcal{M}(B_{2}\otimes A)$ are actions of $G=(A,\Delta)$ on $B_{1}$, respectively $B_{2}$, define the direct sum action of $G$ on $B_{1}\oplus B_{2}$ by \[\delta_{1} \oplus \delta_{2}=(i_{1}\otimes i_{A})\circ \delta_{1} \circ p_{1} +(i_{2}\otimes i_{A})\circ \delta_{2} \circ p_{2}\]
where $i_{k}:B_{k}\hookrightarrow B_{1}\oplus B_{2}$ is the inclusion ($k$=1,2),
             $p_{k}:B_{1}\oplus B_{2}\rightarrow B_{k}$ is the projection ($k$=1,2) 
and  $i_{A}:A\rightarrow A$ is the identity map.  

\end{definition}

\ Notice that $i_{k}\otimes i_{A}$ can be extended as maps from $\mathcal{M}(B_{k}\otimes A)$ to \\ $\mathcal{M}((B_{1}\oplus B_{2})\otimes A)$.

\ If $d: D(d)\rightarrow B$ is a closed, densely defined, $\ast$-derivation, denote by $\Gamma_{d}$ its graph, i.e. :
\[ \Gamma_{d}=\{ x\oplus d(x)\mid x\in D(d)\}\subseteq B\oplus B.\]
\newtheorem{proposition}[definition]{Proposition}
\begin{proposition} d is a $\delta$-invariant derivation if and only if $(\delta \oplus \delta)(\Gamma_{d})\subseteq \overline{\Gamma_{d} \odot A}$, where $\overline{\Gamma_{d} \odot A}$ is the closure of the algebraic tensor product $\Gamma_{d} \odot A$. 
\end{proposition}
\begin{proof} Straightforward from definitions.
\end{proof}

\begin{lemma}\label{L:den} If $d$ is a $\delta$-invariant derivation then $P_{\alpha}(D(d))\subseteq D(d)$ and $P_{\alpha}(D(d))$ are dense in $B_{\alpha}$.
\end{lemma}
\begin{proof} Let $b\in D(d)$. Since $d$ is $\delta$-invariant, 
\[(\delta \oplus \delta)(b\oplus d(b))\in \overline{\Gamma_{d} \odot A}.\]

\ Therefore 
\[ (i\otimes h_{\alpha})((\delta \oplus \delta )(b\oplus d(b))\in \Gamma_{d}\]

\ But this means in particular that $P_{\alpha}(b)\in D(d)$. Then, since $P_{\alpha}$ are continuous maps, we get
\[ B_{\alpha}=P_{\alpha}(B)=P_{\alpha}(\overline{D(d)})\subseteq \overline{P_{\alpha}(D(d))}\subseteq B_{\alpha}.\]
Hence $P_{\alpha}(D(d))$ are dense in $B_{\alpha}$.
\end{proof}

\newtheorem{assumption}[definition]{Assumption}
\begin{assumption}\label{A:as} For the next results we will assume the following:

\underline{Either}

\ a) Both $\delta$ and $h$ are faithful, in which case, the conditional expectation $P_{\iota}: B\rightarrow B^{\delta}$, $P_{\iota}=(\iota \otimes h)\circ \delta$, is faithful.

\underline{or}

\ b) $B$ is a simple C*-algebra.
\end{assumption}

\begin{remark} Either one of a) and b) implies that for every $x\in B$, $\|x\|=\sup_{\tilde{\varphi}}\|\pi_{\tilde{\varphi}}(x)\|$ where $\tilde{\varphi}$ runs through $\{\tilde{\varphi}=\varphi \circ P_{\iota} \mid \varphi$ is a state of $B^{\delta}\}$ and $\pi_{\tilde{\varphi}}$ is the GNS representation of $B$ with respect to the state $\tilde{\varphi}$.
\end{remark}

\begin{remark}\label{R:cl} Assume that \ref{A:as} holds. If $d:\sum_{\alpha \in \widehat{G}}B_{\alpha} \rightarrow B$ is a not necessarily closed, $\ast$-derivation that commutes with $\delta$ (i.e. $ \delta (d(x))=(d\otimes i)(\delta (x))$ for every $x\in \sum B_{\alpha}$) then $d$ is closable and its closure is $\delta$-invariant.
\end{remark}
\begin{proof} Since $d$ commutes with $\delta$, it follows in particular that $d\circ P_{\iota}=P_{\iota} \circ d$. Since $B_{\iota}=B^{\delta}\subset D(d)$, the restriction of $d$ to $B^{\delta}$ is a bounded derivation. By our assumption  \ref{A:as}, the conditions of (\cite{bra}, Proposition 1.4.11) are satisfied and the conclusion follows.
\end{proof}

\ We give next our main results. For a comprehensive theory of generators we send to \cite{brat} and \cite{brat}.

\newtheorem{theorem}[definition]{Theorem}
\begin{theorem}\label{T:gen} Let $(B, G, \delta)$ be a quantum dynamical system satisfying \ref{A:as}. Let $d: \sum B_{\alpha}\rightarrow B$ be a densely defined $\ast$-derivation that commutes with $\delta$ in the sense of \ref{R:cl}. Then $d$ is closable and its closure is the generator of a one parameter group of $\delta$-invariant automorphisms of $B$.
\end{theorem}
\begin{proof} Notice first that if $x\in B_{\alpha}$ we have:
\begin{align} d(x)&=d((i\otimes h_{\alpha})(\delta (x)))=(d\otimes h_{\alpha})(\delta (x)) \notag \\ 
                  &=(i\otimes h_{\alpha})((\overline{d\otimes i})(\delta(x)))=(i\otimes h_{\alpha})(\delta (d(x))). \notag
\end{align}

\ Hence $d(x)\in B_{\alpha}$. By Remark \ref{R:cl} $d$ is closable and its closure, denoted also $d$ is $\delta$-invariant.

\ Since $B_{\alpha}$, $\alpha \in \widehat{G}$, are closed subspaces of $B$ which are $d$-invariant, they consist of analytic elements for $d$, i.e. for every $x\in B_{\alpha}$, the series $\sum \frac{z^{n}}{n!} d^{n}(x)$ is absolutely convergent for every $z\in \mathbb{C}$. Therefore, by Lemma \ref{L:P}, (iv), $d$ has a dense set of analytic elements.

\ We will check next that $\|(1+\alpha d)(x)\| \geqslant \|x\|$, for every $\alpha \in \mathbb{R}$, $x\in D(d)$ and then apply (\cite{brat}, Theorem 3.2.50) to conclude that $d$ is a generator.

\ Let $\varphi$ be a state of $B^{\delta}$ and $\widetilde{\varphi} = \varphi \circ P_{\iota}$ be the corresponding state of $B$. Let $(\pi _{\widetilde{\varphi}}, H_{\widetilde{\varphi}}, \xi_{\widetilde{\varphi}})$ be the associated GNS representation of $B$. Let $\pi = \underset{\varphi}{\oplus} \pi_{\widetilde{\varphi}}$, where $\varphi$ is a state of $B^{\delta}$. Then, by \ref{A:as}, $\pi$ is a faithful representation of $B$ on $H_{\pi}=\underset{\varphi}{\oplus} H_{\pi_{\widetilde{\varphi}}}$. We will identify $B$ with $\pi (B)\subset B(H_{\pi})$. 

\ Since the restriction of $d$ to $B^{\delta}$ is a bounded derivation, by Sakai's Theorem, its extension to the weak closure $\overline{B^{\delta}}^{w}$ is an inner derivation. Let $h_{0}\in \overline{B^{\delta}}^{w}$ be a selfadjoint element such that $d|_{B^{\delta}}= ad(ih_{0})|_{B^{\delta}}$. On the other hand, since, in particular, $B^{\delta}\subset D(d)$, by the von Neumann algebra version of (\cite{bra}, Proposition 1.4.11) it follows that $d$ is $\sigma$-weakly closable on $\overline{B}^{w}$. Let $\overline{d}$ denote the $\sigma$-weak closure of $d$.

\ It is clear that $h_{0}\in D(\overline{d})$ and $\overline{d} (h_{0})=0$. Therefore $\overline{d}\circ ad(ih_{0})=ad(ih_{0})\circ \overline{d}$.

\ Set $d_{1}=\overline{d}-ad(ih_{0})$. Then $d_{1}$ satisfies the conditions of (\cite{bra}, Corollary 1.5.6) and as in the proof of that Corollary, there is a skew symmetric operator $S$ on $\underset{\varphi}{\sum} D(d_{1})\xi_{\widetilde{\varphi}}$ such that:
\begin{align} &1)\quad S(x \xi_{\widetilde{\varphi}}) = d_{1}(x)\xi_{\widetilde{\varphi}},\; \forall \varphi \; \text{state on}\; B^{\delta},\; \forall x \in D(d_{1}).\notag \\
&2)\quad d_{1}(x)=Sx-xS,\; \forall x\in D(d_{1}).\notag
\end{align}

\ Therefore, $\overline{d}=ad(S+ih_{0})$. Since $\underset{\alpha}{\sum}B_{\alpha}$ is a set of analytic elements of $d$, $h_{0}$ is bounded and $Sh_{0}=h_{0}S$, it follows that $\underset{\alpha}{\sum}B_{\alpha}\xi_{\widetilde{\varphi}}$ are analytic elements of $S+ih_{0}$ for every $\varphi$ and therefore the skew symmetric operator $S+ih_{0}$ has a dense set of analytic elements. Indeed, if $x\in B_{\alpha}$ then
\begin{align} \|(S+ih_{0})^{n}x\xi_{\widetilde{\varphi}}\|&=\|\sum_{k=1}^{n} \binom{n}{k}i^{k}(Adh_{0})^{k}S^{n-k}x\xi_{\widetilde{\varphi}}\| \notag \\ 
&= \|\sum_{k=1}^{n} \binom{n}{k}i^{k}(Adh_{0})^{k}d_{1}^{n-k}(x)\xi_{\widetilde{\varphi}}\| \notag \\ 
&\leqslant (2\|h_{0}\|+\rho)^{n}\|x\|  \notag
\end{align}
where $\rho=\|d_{1}|_{B_{\alpha}}\|<\infty$.
Therefore $S+ih_{0}=iK$, where $K$ is essentially selfadjoint.

\ Applying (\cite{brat}, Corollary 3.2.56) it follows that $\|(1+\alpha \overline{d})(x)\| \geqslant \|x\|$ for every $x\in D(\overline{d})$ and the proof is complete.
\end{proof}

\newtheorem{corollary}[definition]{Corollary}
\begin{corollary} If $(B, G, \delta)$ is an ergodic system satisfying Assumption \ref{A:as} and $d: D(d) \rightarrow B$ is a $\delta$-invariant $\ast$-derivation, then $d$ is a generator.
\end{corollary}
\begin{proof} By (\cite{boca}, Theorem ), the spectral subspaces $B_{\alpha}$ are finite dimensional. By Lemma \ref{L:den}, $P_{\alpha}(D(d)) \subseteq B_{\alpha}\cap D(d)$ and $P_{\alpha}(D(d))$ is dense in $B_{\alpha}$. Therefore $P_{\alpha}(D(d))= B_{\alpha}\subset D(d)$ and the conditions of the Theorem \ref{T:gen} are satisfied.
\end{proof}

\ Our next result refers to quantum group actions on the C*-algebra of compact operators $K(H)\subset B(H)$. This result is an extension of (\cite{gdm}, Theorem 4.1). The condition we use is also slightly weaker than the tangential condition used in \cite{gdm}. We do not require that $K(H)^{\delta}\subseteq ker(d)$.

\begin{theorem} Let $\delta$ be an action of $G=(A,\Delta)$ on $K(H)$. If $d:D(d)\rightarrow K(H)$ is a closed densely defined $\ast$-derivation such that:
\begin{align} &a)\; d \; is \; \delta -invariant,\; and\notag \\
&b)\; K(H)^{\delta}\subseteq D(d).\notag
\end{align}
Then $\delta$ is a generator.
\end{theorem}
\begin{proof} We will check the conditions A1), B2) and C1) of \cite{brat}, Theorem 3.2.50.

\ Condition A1) is satisfied by hypothesis. We will show next that $d$ has a dense set of analytic elements (condition B2)).

\ Since $K(H)^{\delta}\subseteq D(d)$ and $d$ is $\delta$-invariant, it follows that $d$ is a (bounded) derivation on $K(H)^{\delta}$. Since every bounded derivation is inner in the $\sigma$-closure, there exists a self-adjoint operator $h_{0} \in \overline{K(H)^{\delta}}^{\sigma}$ such that $d|_{K(H)^{\delta}}=ad(ih_{0})|_{K(H)^{\delta}}$.

\ Notice that, by (\cite{brat}, Corollary 3.2.27), $d$ is $\sigma$-weakly closable on $B(H)$ and $\overline{d}(h_{0})=0$, where $\overline{d}$ denotes the $\sigma$-weak closure of $d$.

\ Consider the derivation $d_{0}=ad(ih_{0})$ on $B(H)$. Then $d_{0}(D(d))\subseteq D(d)$. Indeed, let $x\in D(d)$ and $(e_{\lambda})_{\lambda}$ an approximate identity of $K(H)^{\delta}$. By Lemma \ref{L:app}, $(e_{\lambda})$ is an approximate identity of $K(H)$.

\ Then
\[ x= \underset{\lambda}{norm-lim(e_{\lambda}x)}\] and \[ h_{0}x=\underset{\lambda}{norm-lim(h_{0}e_{\lambda}x)}.\]

\ This last equality follows since $h_{0}e_{\lambda}\in K(H)$, $\sigma - lim h_{0}e_{\lambda}=h_{0}$ and $x\in K(H)$.

\ By Remark \ref{rem}, $\delta(h_{0}e_{\lambda})=\delta''(h_{0})\delta(e_{\lambda})=h_{0}e_{\lambda}\otimes 1$ and hence $h_{0}x\in K(H)^{\delta}$. Therefore $h_{0}x$ (and $xh_{0}$) are elements of $D(d)$, so $d_{0}(D(d))\subseteq D(d)$.

\ Moreover, since $\delta''(h_{0})=h_{0}\otimes 1$ it follows that $d_{0}$ is $\delta$-invariant.

\ It is obvious that $d$ commutes with $d_{0}$. since $d$ is closed and $d_{0}$ bounded, it follows that $d_{1}=d-d_{0}$ is a closed, densely defined (on $D(d)$), $\delta$-invariant derivation. Clearly $d_{1}|_{K(H)^{\delta}}=0$.

\ Let $(e_{\lambda})_{\lambda}$ be an approximate identity of $K(H)^{\delta}$ consisting of projections. Then, by Lemma \ref{L:app} $(e_{\lambda})_{\lambda}$ is an approximate identity of $K(H)$.

\ Since $e_{\lambda}$ are finite dimensional, $e_{\lambda}K(H)e_{\lambda}$ are finite dimensional C*-subalgebras of $K(H)$. Further, since $d_{1}|_{K(H)^{\delta}}=0$, it follows that $d_{1}(e_{\lambda}K(H)e_{\lambda})\subseteq e_{\lambda}K(H)e_{\lambda}$, $ \forall \lambda \in \Lambda$. Therefore $\underset{\lambda \in \Lambda}{\bigcup}e_{\lambda}K(H)e_{\lambda}$ is a dense set of analytic elements of $d_{1}$. Since $d_{0}$ commutes with $d_{1}$ and $d_{0}$ is bounded, it follows that $d$ has a dense set of analytic elements. Indeed, if $x\in e_{\lambda}\mathcal{K}(H)e_{\lambda}$ we have
\begin{align} \|d^{n}(x)\|&=\|(d_{0}+d_{1})^{n}(x)\|\notag\\
&=\|\sum_{k=1}^{n} \binom{n}{k}d_{0}^{k}d_{1}^{n-k}(x)\|
\leqslant (2\|h_{0}\|+\rho)^{n}\|x\| \notag
\end{align}
where $\rho=\|d_{1}|_{e_{\lambda}\mathcal{K}(H)e_{\lambda}}\|<\infty$ (since $e_{\lambda}$ is finite dimensional and $d_{1}(e_{\lambda}\mathcal{K}(H)e_{\lambda})\subseteq e_{\lambda}\mathcal{K}(H)e_{\lambda}$).

\ By (\cite{bra}, Example 1.6.4) and its proof, there is a skew symmetric operator $S$ such that $d_{1}\subseteq ad(S)$ and $S$ has a dense set of analytic elements. Therefore $S$ is essentially skew adjoint.

\ Since $h_{0}$ is a bounded selfadjoint operator, we have that $S+ih_{0}$ is essentially skew adjoint.

\ Since $d\subseteq ad(S+ih_{0})$, by (\cite{bra}, Corollary 1.5.6) it follows that the condition C1) of (\cite{brat}, Theorem 3.2.50) is satisfied and we are done.

\end{proof}

\end{document}